\documentclass{amsart}
\usepackage{url}
\usepackage{enumerate}

\newtheorem{thrm}{Theorem}[section]
\newtheorem{lem}[thrm]{Lemma}

\newtheorem{cor}[thrm]{Corollary}
\theoremstyle{definition}

\newtheorem{remark}[thrm]{Remark}
\numberwithin{equation}{section}
\newcommand{\StirSubset}[2]{ \genfrac{\{}{\}}{0pt}{}{#1}{#2} }

\author{G.~A. Kalugin and D.~J. Jeffrey}
\address{
Department of Applied Mathematics \\
The University of Western Ontario \\
London, Ontario, Canada N6A 5B7}
\email{ gkalugin@uwo.ca, djeffrey@uwo.ca}

\begin{document}

\title[Unimodal sequences]{Unimodal sequences show Lambert W is Bernstein}

\begin{abstract}
We consider a sequence of polynomials appearing in expressions for the derivatives of the Lambert W function.
The coefficients of each polynomial are shown to form a positive sequence that is
log-concave and unimodal. This property implies that the positive real branch of the Lambert W
function is a Bernstein function.
\end{abstract}
\maketitle

\section{Introduction} \label{sect1}
The Lambert $W$ function was defined and studied in~\cite{Corless1996Lam}.
It is a multivalued function having branches $W_k(z)$, each of which obeys
$W_k\exp(W_k) = z$.
The principal branch $W_0$ maps the set of positive reals to itself, and is the only
branch considered here. Therefore we omit the subscript $0$ for brevity.
The $n$th derivative of $W$ is given implicitly by
\begin{equation}\label{eq:nthderiv}
 \frac{d^n W(x)}{d x^n} =  \frac{ \exp(-n W(x)) p_n(W(x)) } { (1 + W(x))^{2n-1} }
\mbox{\qquad for\ } n\ge 1\ ,
\end{equation}
where the polynomials $p_n(w)$ satisfy $p_1(w)=1$, and the recurrence relation
\begin{equation}
\label{eq:polyrecurr}
   p_{n+1}(w) = -(nw + 3n - 1)p_n(w) + (1+w)p_n'(w)\mbox{ \qquad for\ }n\ge 1\ .
\end{equation}
In \cite{Corless1997seqseries}, the first 5 polynomials were printed explicitly:
\begin{align*}
p_1(w)&=1\ ,\ p_2(w)=-2-w\ ,\ p_3(w) = 9+8w+2w^2\ ,\\
p_4(w)&=-64-79w-36w^2-6w^3\ ,\\
p_5(w)&=625 + 974w+ 622w^2+ 192w^3+ 24w^4\ .
\end{align*}
These initial cases suggest the conjecture that each polynomial $(-1)^{n-1}p_n(w)$ has all positive coefficients,
and if this is true, then $dW(x)/dx$ is a completely monotonic function~\cite{SokalComm2008}.
We here prove the conjecture and prove in addition that the coefficients are unimodal and log-concave.

\section{Formulae for the coefficients}\label{sect2}
In view of the conjecture, we write
\begin{equation}\label{eq:Polydef}
   p_n(w) = (-1)^{n-1}\sum_{k=0}^{n-1}\beta_{n,k}w^k\ .
\end{equation}
We now give several theorems regarding the coefficients.
\begin{thrm}\label{th:Formulae}
The coefficients $\beta_{n,k}$ defined in \eqref{eq:Polydef} obey the recurrence relations
\begin{align}
 &\beta_{n,0}=n^{n-1}\ ,\quad \beta_{n,1}=3n^n-(n+1)^n-n^{n-1} \ , \label{eq:begin} \\
 &\beta_{n,n-1}=(n-1)!\ ,\quad \beta_{n,n-2}=(2n-2)(n-1)! \ , \label{eq:end} \\
 &\beta_{n+1,k}=(3n-k-1) \beta_{n,k} +n\beta_{n,k-1}-(k+1)\beta_{n,k+1}\ , \quad 2\le k \le n-3 \ . \label{eq:recurrence}
\end{align}
\end{thrm}
\begin{proof}
By substituting (\ref{eq:Polydef}) into (\ref{eq:polyrecurr}) and equating coefficients.
\end{proof}

\begin{thrm}
An explicit expression for the coefficients $\beta_{n,k}$ is
\begin{equation}\label{eq:ExplicitOne}
\beta_{n,k}= \sum_{m=0}^k\frac{1}{m!} \binom{2n-1}{k-m} \sum_{q=0}^m \binom{m}{q} (-1)^q (q+n)^{m+n-1} \ .
\end{equation}
\end{thrm}
\begin{proof}
We rewrite (\ref{eq:nthderiv}) in the form
\[ p_n(W(x)) = (1+W(x))^{2n-1} e^{nW(x)} \frac{d^n W(x)}{dx^n}\ .
\]
From the Taylor series of $W(x)$ around $x=0$, given in \cite{Corless1996Lam}, we obtain
\[ \frac{d^n W(x)}{dx^n} = \sum_{m=n}^\infty \frac{(-m)^{m-1}}{(m-n)!} x^{m-n}\ .
\]
Substituting this into the expression of $p_n$, using $x=We^W$ and changing the index of summation, we obtain the equation
\begin{equation}\label{eq:polynom}
p_n(w) = (1+w)^{2n-1}\sum_{s=0}^\infty(-1)^{n+s-1}(n+s)^{n+s-1}\frac{w^s}{s!}e^{(n+s)w}\ .
\end{equation}
We expand the right side around $w=0$ and equate coefficients of $w$.
\end{proof}

\begin{remark}\label{r1}
The polynomials $p_n(w)$ can be expressed in terms of the diagonal Poisson transform
$\mathbf{D}_n[f_s;z]$ defined in \cite{PobleteViolaMunro1997}, namely, by \eqref{eq:polynom}
\begin{equation} \label{eq:DPT}
p_n(w)=(-1)^{n-1}(1+w)^{2(n-1)}\mathbf{D}_n[(n+s)^{n-1};-w]\ .
\end{equation}
\end{remark}



\begin{thrm} \label{Th:specnum}
The coefficients can equivalently be expressed either in terms of shifted
$r$-Stirling numbers of the second kind $\StirSubset{n+r}{m+r}_r$ defined in \cite{Broder1984},
\begin{equation}
\label{eq:r-Stir}
\beta_{n,k} = \sum_{m=0}^k(-1)^m \binom{2n-1}{k-m}
\StirSubset{2n -1 + m}{ n + m}_n\quad \ ,
\end{equation}
or in terms of Bernoulli polynomials of higher order $B_n^{(z)}(\lambda)$ defined in \cite{Norlund1924},
\begin{equation} \label{eq:Bernoulli}
\beta_{n,k}= \sum_{m=0}^k(-1)^m \binom{2n-1}{k-m}
\binom{m+n-1}{n-1}B_{n-1}^{(-m)}(n) \ ,
\end{equation}
or in terms of the forward difference operator $\Delta$ \cite[p.\:188]{GrahamKnuthPatashnik},
\[
\beta_{n,k}=\sum_{m=0}^k
\binom{2n-1}{k-m}
\frac{(-1)^m}{m!}\Delta^m n^{m+n-1} \ .
\]

\end{thrm}
\begin{proof}
We convert \eqref{eq:ExplicitOne} using identities found in \cite{Broder1984} and \cite{LopezTemme2010} respectively.
\[
\StirSubset{n+r}{m + r}_r= \frac{1}{m!}\sum_{q=0}^m (-1)^{m-q}\binom{m}{q}(q+r)^n
\]
and
\[
B_n^{(-m)}(r)= \frac{n!}{(m+n)!}\sum_{q=0}^m (-1)^{m-q}\binom{m}{q}(q+r)^{m+n}\ .
\]
\end{proof}

\section{Properties of the coefficients}\label{sect3}
We now give theorems regarding the properties of the $\beta_{n,k}$.
We recall the following definitions \cite{Stanley1989}.
A sequence $c_0, c_1 \ldots c_n$ of real numbers is said to be \emph{unimodal} if for some $0\le j\le n$  we have
$c_0\le c_1\le \ldots \le c_j \ge c_{j+1} \ge \ldots \ge c_n$, and it is said to be \emph{logarithmically concave}
(or log-concave for short) if $c_{k-1}c_{k+1} \le c_k^2$ for all $1\le k\le n-1$.
We prove that for each fixed $n$, the $\beta_{n,k}$ are unimodal and log-concave with respect to $k$.
Since a log-concave sequence of positive terms is unimodal \cite{Wilf2005}, it is convenient to start with the log-concavity property.

\begin{thrm} \label{Th:LC}
For fixed $n\geq3$ the sequence $\left\{k!\beta_{n,k}\right\}_{k=0}^{n-1}$ is log-concave.
\end{thrm}
\begin{proof}
Using \eqref{eq:ExplicitOne} we can write
\begin{equation}
k!\beta_{n,k}=(2n-1)!\sum_{m=0}^k \binom{k}{m}x_m y_{k-m} \ , \notag
\end{equation}
where
\begin{equation}\label{eq:x}
x_m= \sum_{j=0}^j \binom{m}{j}a_j \ , \quad a_j=(-1)^j (n+j)^{m+n-1} \ ,
\end{equation}
and $y_m=1/(2n-1-m)!$ .
Since the binomial convolution preserves the log-concavity property \cite{Walkup1976,WangYeh2007LC},
it is sufficient to show that the sequences $\left\{x_m\right\}$ and $\left\{y_m\right\}$ are log-concave.
We have
\begin{equation}
 \begin{split}
  a_{j-1}a_{j+1}=(-1)^{j-1}(n+j-1)^{m+n-1}(-1)^{j+1}(n+j+1)^{m+n-1} \\
  =(-1)^{2j}\left((n+j)^2-1\right)^{m+n-1}<(-1)^{2j}(n+j)^{2(m+n-1)}=a_j^2 \ . \notag
 \end{split}
\end{equation}
Thus the sequence $\left\{a_j\right\}$ is log-concave and so is $\left\{x_m\right\}$
due to \eqref{eq:x} and the afore-mentioned property of the binomial convolution.
The sequence $\left\{y_m\right\}$ is log-concave because
\[
 \begin{split}
  y_{m-1}y_{m+1}&=\frac{1}{(2n-1-m+1)!}\frac{1}{(2n-1-m-1)!} \\
                &=\frac{2n-1-m}{2n-1-m+1}\frac{1}{(2n-1-m)!}\frac{1}{(2n-1-m)!}<y_m^2 \ .
  \end{split}
\]
\end{proof}
Now we prove that the coefficients $\beta_{n,k}$ are positive.
The following two lemmas are useful.
\begin{lem} \label{Lemma:Gen}
If a positive sequence $\left\{k!c_k\right\}_{k\geq0}$ is log-concave, then
\renewcommand{\labelenumi}{\upshape{(\roman{enumi})}}
\begin{enumerate}
\item $\left\{(k+1)c_{k+1}/c_k\right\}$ is non-increasing;
\item $\left\{c_k\right\}$ is log-concave;
\item the terms $c_k$ satisfy
\end{enumerate}
\begin{equation} \label{ineq:genlemma}
c_k c_m \geq \binom{k+m}{k} c_0 c_{k+m} \quad (0 \leq m \leq k+1) \ .
\end{equation}
\end{lem}
\begin{proof}
The statements (i) and (ii) are obvious. To prove (iii) we apply a method used in \cite{AsaiKuboKuo2000}.
Specifically, by (i) we have for $0\leq p\leq k$
\begin{equation}
\frac{c_{p+1}}{c_p}\geq\frac{k+p+1}{p+1}\frac{c_{k+p+1}}{c_k} \ . \notag
\end{equation}
Apply the last inequality for $p=0,1,2,...,m$ with $m \leq k+1$, and
form the products of all left-hand and right-hand sides.
As a result, after the cancellation we obtain
\begin{equation*}
\frac{c_m}{c_0}\geq\frac{k+1}{1}\frac{k+2}{2} \ldots \frac{k+m}{m}\frac{c_{k+m}}{c_k} \ ,
\end{equation*}
which is equivalent to \eqref{ineq:genlemma}.
\end{proof}
\begin{lem}
If the coefficients $\beta_{n,k}$ are positive, then for fixed $n\geq3$ they satisfy
\begin{equation} \label{ineq:lemma}
\frac{(k+1)\beta_{n,k+1}}{\beta_{n,k}}<n-1 \ .
\end{equation}
\end{lem}
\begin{proof}
By Theorem \ref{Th:LC} and under the assumption of lemma, for fixed $n\geq3$
the sequence $\left\{k!\beta_{n,k}\right\}_{k=0}^{n-1}$ meet the conditions of Lemma \ref{Lemma:Gen}.
Applying the inequality \eqref{ineq:genlemma} with $m=1$ to this sequence gives $(k+1)\beta_{n,k+1}/\beta_{n,k}\leq \beta_{n,1}/\beta_{n,0}$.
Then the lemma follows as due to \eqref{eq:begin}
\begin{equation}
\frac{\beta_{n,1}}{\beta_{n,0}}=\frac{3n^n-(n+1)^n-n^{n-1}}{n^{n-1}}=3n-n\left(1+\frac{1}{n}\right)^n-1<3n-2n-1=n-1 \ . \notag
\end{equation}
\end{proof}
\begin{thrm} \label{Th:positive}
The coefficients $\beta_{n,k}$ are positive.
\end{thrm}
\begin{proof}
We prove the statement by induction on $n$.
It is true for $n\le 5$ (see \S 1). Assume that for some fixed $n$
all the members of the sequence $\left\{\beta_{n,k}\right\}_{k=0}^{n-1}$ are positive.
Since $\beta_{n+1,0}=(n+1)^n>0$ and $\beta_{n+1,n}=n!>0$ by \eqref{eq:begin} and \eqref{eq:end},
we only need to consider $k=1,2,...,n-1$.

Substituting inequalities $\beta_{n,k+1}<(n-1)\beta_{n,k}/(k+1)$ and $\beta_{n,k-1}>k\beta_{n,k}/(n-1)$,
which follow from \eqref{ineq:lemma}, in the recurrence \eqref{eq:recurrence} immediately gives the result
\begin{equation*}
\beta_{n+1,k}>(3n-k-1)\beta_{n,k}+n\frac{k}{n-1}\beta_{n,k}-(k+1)\frac{n-1}{k+1}\beta_{n,k}
   =\left(2n+\frac{k}{n-1}\right)\beta_{n,k}>0 \ .
\end{equation*}
Thus the proof by induction is complete.
\end{proof}
\begin{cor}
The sequence $\left\{\beta_{n,k}\right\}_{k=0}^{n-1}$ is unimodal for $n\geq3$.
\end{cor}
\begin{proof}
By Theorem \ref{Th:positive} the sequence $\left\{\beta_{n,k}\right\}_{k=0}^{n-1}$ is positive,
therefore by Theorem \ref{Th:LC} and Lemma \ref{Lemma:Gen}(ii) it is log-concave and, hence, unimodal.
\end{proof}

\section{Relation to Carlitz's numbers}\label{sect4}
There is a relation between the coefficients $\beta_{n,k}$ and numbers $B(\kappa,j,\lambda)$ introduced
by Carlitz in \cite{Carlitz1980WII}. Comparing the formulae \eqref{eq:r-Stir} and \eqref{eq:inverse}
with the corresponding \cite[eq.(6.3)]{Carlitz1980WII} and \cite[eq.(2.9)]{Carlitz1980WII},
taking into account that he uses the notation $R(n,m,r)=\StirSubset{n+r}{m+r}_r$, we find
\begin{equation} \label{eq:termsB}
\beta_{n,k}=(-1)^kB(n-1,n-1-k,n) \ .
\end{equation}
It follows that for $n\geq3$, the sequence $\left\{B(n-1,k,n)\right\}_{k=0}^{n-1}$ is log-concave
together with $\left\{\beta_{n,k}\right\}_{k=0}^{n-1}$.

Using the property \cite[eq.(2.7)]{Carlitz1980WII} that $\sum_{j=0}^\kappa B(\kappa,j,\lambda)=(2\kappa-1)!!$,
we can compute $p_n(w)$ at the singular point where $W=-1$ (cf. \eqref{eq:nthderiv}).
Thus, substituting $w=-1$ in \eqref{eq:Polydef} gives $p_n(-1)=(-1)^{n-1}(2n-3)!!$.
Thus $w=-1$ is not a zero of $p_n(w)$.

We also note that the numbers $B(\kappa,j,\lambda)$ are polynomials of $\lambda$
and satisfy a three-term recurrence \cite[eq.\:(2.4)]{Carlitz1980WII}
\begin{equation} \label{eq:recB}
B(\kappa,j,\lambda)=(\kappa+j-\lambda)B(\kappa-1,j,\lambda)+(\kappa-j+\lambda)B(\kappa-1,j-1,\lambda)
\end{equation}
with $B(\kappa,0,\lambda)=(1-\lambda)^{\bar{k}}$,$\quad B(0,j,\lambda)=\delta_{j,o}$.
This gives one more way to compute the coefficients $\beta_{n,k}$, specifically,
for given $n$ and $k$ we find a polynomial \mbox{$B(n-1,n-1-k,\lambda)$}
using \eqref{eq:recB} and then set $\lambda=n$ to use \eqref{eq:termsB}.

\section{Concluding remarks}\label{sect5}
It has been established that the coefficients of the polynomials $(-1)^{n-1}p_n(w)$ are positive, unimodal and log-concave.
These properties imply an important property of $W$.
In particular, it follows from formula \eqref{eq:nthderiv} and Theorem \ref{Th:positive}
that $(-1)^{n-1}(dW/dx)^{(n-1)}>0$ for $n\geq1$.
Since $W(x)$ is positive for all positive $x$~\cite{Corless1996Lam}, this means that the derivative $W'$ is completely
monotonic and $W$ itself is a Bernstein function \cite{Berg2008}.

Some additional identities can be obtained from the results above.
For example, computing $\beta_{n,n-1}$ by \eqref{eq:r-Stir} and comparing with \eqref{eq:end} gives
\[
\sum_{m=0}^{n-1}(-1)^m \binom{2n-1}{n-m-1}
\StirSubset{2n -1 + m}{ n + m}_n
=(n-1)! \ .
\]
A relation between $\StirSubset{2n -1 + m}{ n + m}_n$ and $B_{n-1}^{(-m)}(n)$ can be obtained from
\eqref{eq:r-Stir} and \eqref{eq:Bernoulli}, but this is a special case of \cite[eq.\:(7.5)]{Carlitz1980WII}.
It is finally interesting to note that \eqref{eq:r-Stir} and \eqref{eq:Bernoulli} can be inverted.
Indeed, in these formulae for fixed $n$,
the sequence $(-1)^k\beta_{n,k}$ is a convolution of two sequences and
the corresponding relation between their generating functions is $G(w)=(1-w)^{2n-1}F(w)$.
Since $F(w)=G(w)(1-w)^{-(2n-1)}=G(w)\sum_{k\geq0}\binom{2n-2+k}{2n-2}w^k$,
the inverse of, for example, \eqref{eq:r-Stir} is
\begin{equation} \label{eq:inverse}
\StirSubset{2n -1 + m}{ n + m}_n = \sum_{k=0}^{n-1}(-1)^k\beta_{n,k}\binom{2n-2+m-k}{2n-2} \ .
\end{equation}

\proof[Acknowledgements]
We thank Prof. Alan Sokal for sending us his conjecture and for his interest and encouragement.
The work was supported by a Discovery Grant from the Natural Sciences and Engineering Research Council of Canada.

\end{document}